\documentclass[12pt]{article}

\usepackage{latexsym, amssymb, amsmath, amsfonts, epsfig, graphicx, colordvi,verbatim,ifpdf}
\usepackage{amsfonts, amsmath, amssymb}
\usepackage{amssymb,amsfonts, cite, psfrag, eepic,color}
\usepackage{amscd,graphics}

\usepackage{graphicx}
\usepackage{color}
\usepackage{ifpdf}
\usepackage{fancybox}

\newtheorem{thm}{Theorem}[section]

\newtheorem{defi}[thm]{Definition}
\newtheorem{lem}[thm]{Lemma}

\def\pf{\noindent{\it Proof.} }
\def\qed{\nopagebreak\hfill{\rule{4pt}{7pt}}
\medbreak}

\numberwithin{equation}{section}

\def\qed{\nopagebreak\hfill{\rule{4pt}{7pt}}
\medbreak}

\newlength{\boxedparwidth}
  {\begin{center} \begin{tabular}{|@{\hspace{.315in}}c@{\hspace{.15in}}|}
                  \hline \\ \begin{minipage}[t]{\boxedparwidth}
                  \setlength{\parindent}{.25in}}%
  {\end{minipage} \\ \\ \hline \end{tabular} \end{center}}

\makeatletter
\@addtoreset{figure}{section}
\@addtoreset{table}{section}
\makeatother

\begin{document}
\begin{center}

{\Large \bf The spt-Crank for Ordinary Partitions}
\end{center}

\vskip 5mm
\begin{center}
{  William Y.C. Chen}$^{a,b}$, {Kathy Q. Ji}$^{a}$ and
  {Wenston J.T. Zang}$^{a}$ \vskip 2mm

   $^{a}$Center for Combinatorics, LPMC-TJKLC\\
   Nankai University, Tianjin 300071, P.R. China\\[6pt]
   $^{b}$Center for Applied Mathematics\\
Tianjin University,  Tianjin 300072, P. R. China\\[6pt]

   \vskip 2mm

   Email addresses: chen@nankai.edu.cn, ji@nankai.edu.cn, wenston@mail.nankai.edu.cn
\end{center}

Dedicated to Professor George E. Andrews on the Occasion of His 75th  Birthday

\vskip 6mm \noindent {\bf Abstract.}
 The spt-function  $spt(n)$ was  introduced by Andrews as the weighted counting of partitions of $n$ with respect to the number of occurrences of the smallest
part.  Andrews, Garvan and Liang defined the spt-crank of an $S$-partition which leads to combinatorial interpretations of the congruences of $spt(n)$ mod $5$ and $7$.  Let $N_S(m,n)$ denote the net number of $S$-partitions of $n$ with spt-crank $m$.  Andrews, Garvan and Liang showed that $N_S(m,n)$ is nonnegative for all integers $m$ and positive integers $n$, and they asked the question of
finding a combinatorial interpretation of $N_S(m,n)$. In this paper, we introduce the structure of doubly marked partitions and define the spt-crank of a doubly marked partition. We show that  $N_S(m,n)$ can be interpreted as the number of
doubly marked partitions of $n$ with spt-crank $m$. Moreover,  we
establish a bijection between marked partitions of $n$ and  doubly marked
partitions of $n$. A marked partition is defined by  Andrews, Dyson and Rhoades
as a partition with exactly one of the smallest parts marked.
They consider it a challenge  to find a definition of the
spt-crank of a marked partition so that the set of
marked partitions of $5n+4$ and $7n+5$ can be divided into five and seven equinumerous classes.
The definition of spt-crank for doubly marked partitions and the bijection between the marked partitions and doubly marked partitions leads to a solution to the problem of Andrews, Dyson and Rhoades.



\section{Introduction}

Andrews \cite{Andrews-2008} introduced the spt-function $spt(n)$ as the
weighted counting of partitions with respect to the
number of occurrences of the smallest part and he discovered
that the spt-function bears striking resemblance to the classical
partition function $p(n)$. Much attention has been drawn to the
investigation of the spt-function, in particular, the spt-crank
of an $S$-partition, see, for example, Andrews, Dyson and Rhoades \cite{Andrews-Dyson-Rhoades-2013}, Andrews, Garvan and Liang \cite{Andrews-Garvan-Liang-2011, Andrews-Garvan-Liang}, Folsom and Ono \cite{Folsom-Ono-2008}, Garvan \cite{Garvan} and Ono \cite{Ono-2011}.

 In this paper, we introduce the structure of doubly marked partitions
 and define the spt-crank of a doubly marked partition.
 This gives a solution to a problem posed by Andrews, Garvan and Liang \cite{Andrews-Garvan-Liang-2011}
on the spt-crank of an $S$-partition.
Moreover,  we find a bijection
between marked partitions and doubly marked partitions,
which leads to  a solution to
a problem of finding the definition of
the spt-crank for ordinary partitions  posed by Andrews, Dyson and Rhoades \cite{Andrews-Dyson-Rhoades-2013}.

Let us give an overview of notation and known results on the spt-function $spt(n)$. For a partition $\lambda$ of $n$, we use $n_s(\lambda)$ to denote
the number of occurrences of the smallest part in $\lambda$. Let $P(n)$ denote the set of ordinary partitions of $n$, then we have
\begin{equation}\label{defi-spt-1}
spt(n)=\sum_{\lambda \in P(n)}n_s(\lambda).
\end{equation}
For example, for $n=4$, we have $spt(4)=10$. Partitions in $
P(4)$ and the values of $n_s(\lambda)$ are
listed below:
\[\begin{array}{c|c}
\lambda & n_s(\lambda)\\[2pt] \hline
(4) & 1\\[2pt]
(3,1)& 1\\[2pt]
(2,2)&2\\[2pt]
(2,1,1)&2\\[2pt]
(1,1,1,1)& 4

\end{array}
\]
The spt-function $spt(n)$ can also be interpreted by marked partitions,
   see Andrews, Dyson and Rhoades \cite{Andrews-Dyson-Rhoades-2013}. A marked partition of $n$  means   a pair $(\lambda,k)$ where $\lambda$ is an ordinary  partition of $n$ and $k$ is an integer identifying one of its smallest parts. If there are $s$ smallest parts in $\lambda$, then $k=1,2,\ldots, s$. However, for the purpose of this paper, we shall mark the
   unique smallest part by its index. More precisely,
   if $\lambda_k$ is the marked smallest part, then we use $(\lambda,k)$
   to denote this marked partition.   For example, there are ten marked partitions of $4$.
  \[\begin{array}{lllll}
((4),1),&  ((3,1), 2),& ((2,2),1),&((2,2),2),&((2,1,1),2)\\[5pt]
((2,1,1),3),&((1,1,1,1),1),&((1,1,1,1),2),&((1,1,1,1),3),&((1,1,1,1),4).
\end{array}
\]

 From the definition \eqref{defi-spt-1} of  spt-function, one can derive the following
 generating function
\begin{equation}\label{gf-spt}
\sum_{n\geq 1}spt(n)q^n  =\sum_{n=1}^\infty \frac{q^n}{(1-q^n)^2(q^{n+1};q)_\infty},
\end{equation}
 see Andrews \cite{Andrews-2008}. Based on the above formula,
 Andrews, Garvan and Liang \cite{Andrews-Garvan-Liang-2011}
 noticed that the generating
 function of $spt(n)$ can expressed in the following form
 \begin{equation}\label{gf-spt-1}\sum_{n\geq 0}spt(n)q^n=\sum_{n=1}^{+\infty}\frac{q^n(q^{n+1};q)_\infty}{(q^n;q)_\infty (q^n;q)_\infty},
\end{equation}
 and they introduced the structure of $S$-partitions and interpreted the right-hand side of \eqref{gf-spt-1}  as the
 generating function of the net number of $S$-partitions of $n$,
 that is, the sum of   signs of $S$-partitions of $n$. In other words,
 Andrews, Garvan and Liang established the following relation
 \[spt(n)=\sum_{\pi }\omega(\pi),\]
 where $\pi$ ranges over $S$-partitions of $n$ and $\omega(\pi)$ is the
 sign of $\pi$.

 To be precise, let $\mathcal{D}$ denote  the set of partitions into distinct parts and $\mathcal{P}$ denote the set of partitions. For   $\lambda \in \mathcal{P}$, we use $s(\lambda)$ to denote  the smallest part  of $\lambda$ with the convention that $s(\emptyset)=+\infty$. Let $\ell(\lambda)$ denote the number of parts of $\lambda$.   The set of $S$-partitions is defined by
\[S=\{(\pi_1, \pi_2,\pi_3)\in \mathcal{D} \times \mathcal{P} \times \mathcal{P}\ | \  \pi_1 \neq \emptyset \text{ and }s(\pi_1)\leq \min\{s(\pi_2),s(\pi_3)\}\}.\]
For $\pi=(\pi_1,\,\pi_2,\,\pi_3) \in S$,
 we  define the weight of $\pi$ to be $|\pi_1|+|\pi_2|+|\pi_3|$ and  we associate with $\pi$ a sign \[\omega(\pi)=(-1)^{\ell(\pi_1)-1}.\]

Using the  generating function \eqref{gf-spt} and Watson's $q$-analog of Whipple's theorem \cite[p.43, eq. (2.5.1)]{Gasper-Rahman-2004},
 Andrews showed that the spt-function can be expressed in terms of the
second moment $N_2(n)$  of ranks, namely,
\begin{align}\label{spt-moments}
spt(n)&=np(n)-\frac{1}{2}N_2(n).
\end{align}
 In general, the $k$th moment  $N_k(n)$ of ranks
 was introduced by Atkin and Garvan \cite{Atkin-Garvan-2011} as  given by
 \begin{eqnarray*}
N_k(n)&=& \sum_{m=-\infty}^{+\infty}m^kN(m,n),
\end{eqnarray*}
where $N(m,n)$  is the number of partitions of $n$ with rank $m$, and the
rank of a partition  is defined as the largest part minus the
number of parts.

In view of relation \eqref{spt-moments} and identities
on the refinements of $N(m,n)$, Andrews proved that
 $spt(n)$ satisfies  congruences mod $5$, $7$ and $13$
 reminiscent to Ramanujan's congruences for $p(n)$. To be more specific, let $N(i,t,n)$ denote the number of partitions of $n$ with rank congruent $i\mod t$. Dyson \cite{Dyson-1944} conjectured
    \begin{eqnarray}\label{Dyson-1}
    N(i,5,5n+4)&=&\frac{p(5n+4)}{5}\quad\text{for } \quad 0\leq i\leq 4,\\[3pt]
    N(i,7,7n+5)&=&\frac{p(7n+5)}{7}\quad\text{for } \quad 0\leq i\leq 6. \label{Dyson-2}
    \end{eqnarray}
These relations were confirmed by Atkin and Swinnerton-Dyer \cite{ Atkin-Swinnerton-Dyer-1954} which imply Ramanujan's congruences mod $5$ and $7$ for $p(n)$ \cite{Ramanujan-1919}:
\begin{eqnarray}
p(5n+4)&\equiv& 0 \pmod 5, \label{pcon-5}\\[3pt]
p(7n+5) &\equiv& 0 \pmod 7,\label{pcon-7} \\[3pt]
p(11n+6) &\equiv& 0 \pmod {11}.  \label{pcon-11}
\end{eqnarray}
Using relation \eqref{spt-moments} along with \eqref{Dyson-1} and \eqref{Dyson-2}, Andrews \cite{Andrews-2008} showed that
    \begin{eqnarray}
 spt(5n+4)&\equiv& 0 \pmod 5, \label{con-5}\\[3pt]
spt(7n+5) &\equiv& 0 \pmod 7 \label{con-7}.
\end{eqnarray}
Andrews also proved the following congruence
\begin{equation}\label{con-13}
spt(13n+6) \equiv 0 \pmod {13}
\end{equation}
by using identities on  $N(i,13,13n+6)$  due to O'Brien   \cite{O'Brien-1965}.
Let
\[r_{a,b}(d)=\sum_{n=0}^\infty(N(a,13,13n+d)-N(b,13,13n+d))q^{13n},\]
and for $1\leq i \leq 5$, let
\[
S_i(d)=r_{(i-1),i}(d)-(7-i)r_{5,6}(d).
\]
O' Brien proved that
\begin{equation}\label{identity-mod13-1}
S_1(6)+2S_2(6)-5S_5(6)\equiv 0\pmod{13}
\end{equation}
and
\begin{equation}\label{identity-mod13-2}
S_2(6)+5S_3(6)+3S_4(6)+3S_5(6)\equiv 0\pmod{13}.
\end{equation}
By relation \eqref{spt-moments}, Andrews derived an expression
of $spt(13n+6)$ in terms of the numbers $N(i,13, 13n+6)$ mod 13. Then the congruence \eqref{con-13} follows from  \eqref{identity-mod13-1} and \eqref{identity-mod13-2}.

 To give combinatorial interpretations of the spt-congruences \eqref{con-5} and \eqref{con-7},  Andrews, Garvan and Liang \cite{Andrews-Garvan-Liang-2011} defined the   spt-crank for $S$-partitions, which  takes the same
 form as the crank for vector partitions. Recall that
 the crank for vector partitions has been used to interpret  Ramanujan's congruences  for $p(n)$ mod $5$, $7$ and $11$,  see Andrews and Garvan \cite{Andrews-Garvan-1988}, Dyson \cite{Dyson-1989} and Garvan \cite{Garvan-1988}.
 Let $\pi$ be an $S$-partition, the spt-crank of $\pi$, denoted  $r(\pi)$, is defined to be the number of parts of $\pi_2$ minus the number of parts of
  $\pi_3$, that is,
 \[ r(\pi)=\ell(\pi_2)-\ell(\pi_3).\]
 Let $N_S(m,n)$ denote the net number of $S$-partitions of $n$ with spt-crank $m$, that is,
\begin{equation}
N_S(m,n)=\sum_{\stackrel{|\pi|=n}{r(\pi)=m}}\omega(\pi)
\end{equation}
and let $N_S(k,t,n)$  denote the net number of $S$-partitions of $n$ with spt-crank congruent $k \pmod t$, namely,
\[N_S(k,t,n)=\sum_{m \equiv k \pmod{t}}N_S(m,n).\]

Andrews, Garvan and Liang \cite{Andrews-Garvan-Liang-2011}  established the following relations.

\begin{thm}\label{AGL} For $0\leq k\leq 4$, we have
\begin{eqnarray*}
N_S(k,5,5n+4)&=& \frac{spt(5n+4)}{5},
\end{eqnarray*}
and for $0\leq k\leq 6$, we have
\begin{eqnarray*}
N_S(k,7,7n+5)&=& \frac{spt(7n+5)}{7}.
\end{eqnarray*}
\end{thm}

By using generating functions, Andrews, Garvan and Liang \cite{Andrews-Garvan-Liang-2011} obtained the
following positivity result for $N_S(m,n)$.

\begin{thm}For all integers $m$ and positive integers $n$, we have
\begin{equation}
N_S(m,n)\geq 0.
\end{equation}
\end{thm}

Dyson \cite{Dyson-2012} gave an alternative proof of this fact by using the following recurrence relation.
\[N_S(m,n)=\sum_{k=1}^\infty (-1)^{k-1}\sum_{j=0}^{k-1}p(n-k(m+j)-(k(k+1)/2)).
\]

Andrews, Garvan and Liang \cite{Andrews-Garvan-Liang-2011} asked the question of finding a combinatorial interpretation of $N_S(m,n)$.
Using the generating function for $N_S(m,n)$ given by Andrews, Garvan and Liang, we  show that $N_S(m,n)$ can be  interpreted as the number of doubly marked partitions of $n$ with spt-crank $m$. This gives a solution to the problem of Andrews, Garvan and Liang.

 Andrews, Dyson and Rhoades \cite{Andrews-Dyson-Rhoades-2013}  proposed the problem of  finding a definition of the spt-crank for marked partitions so that the set of marked partitions of $5n+4$ and $7n+5$ can be divided into five and seven equinumerous classes.
We establish a bijection $\Delta$ between the set of marked partitions of $n$ and the set of doubly marked partitions of $n$. While the spt-crank of a doubly
marked partition does not lead to an explicit formula
in terms of the corresponding marked partition, there is
a way to define the spt-crank of a marked partition based on the
bijection $\Delta$ between marked partitions and doubly marked partitions. Hence, in principle,
the spt-crank of a doubly marked partition can
be considered as a solution to the problem of Andrews, Dyson and Rhoades.
It would be interesting to find an spt-crank that can be directly
defined on marked partitions.

\section{Combinatorial interpretation of $N_S(m,n)$}

In this section, we first define the doubly marked partitions and the spt-crank of a doubly marked partition. Then we show that  $N_S(m,n)$ equals the number of  doubly marked partitions of $n$ with spt-crank $m$. In the following
definition, we assume that a partition $\lambda$ of $n$ is represented by
its Ferrers diagram, and we use $D(\lambda)$ to denote size of the
Durfee square of $\lambda$, see \cite[p. 28]{Andrews-1976}.

 We now give the definition of doubly marked partitions.

\begin{defi}
A doubly marked partition of $n$ is an ordinary partition $\lambda$ of $n$ along with two distinguished columns indexed by $s$ and $t$, denoted $(\lambda,s,t)$, where
\begin{itemize}
\item[{\rm (1)}] $1\leq s\leq D(\lambda)${\rm;}
\item[{\rm (2)}] $s\leq t\leq \lambda_1${\rm;}
\item[{\rm (3)}] $\lambda'_s=\lambda'_t$.
\end{itemize}
\end{defi}

For example, $((3,2,2),1,2)$ is a doubly marked partition, whereas $((3,2,1),1,2)$ and $((3,2,2),2,1)$ are not doubly marked partitions, see Figure 2.1.
\input{double.TpX}

 To define the spt-crank of a doubly marked partition  $(\lambda,s,t)$,
 let \begin{equation}\label{defi-t}
  g(\lambda,s,t)=\lambda'_{s}-s+1,
\end{equation}
where $\lambda'$ denotes the conjugate of $\lambda$, in other words,
 $\lambda'_s$ is the number of parts in $\lambda$ that are not less than $s$.   Since $s\leq D(\lambda)$, we see that $\lambda'_s\geq s$, which implies that $g(\lambda,s,t)\geq 1$.

\begin{defi}Let $(\lambda,s,t)$ be a doubly marked partition, and
 let $g=g(\lambda, s,t)$. The  spt-crank of   $(\lambda,s,t)$ is defined by
\begin{equation}\label{def-sptc}
c(\lambda,s,t)=g-\lambda_{g}+t-s.
\end{equation}
\end{defi}

For example, for the doubly marked partition $((4,4,1,1),2,3)$, we have $g=2-1=1$ and the spt-crank  equals $1-\lambda_1+3-2=-2.$

The following theorem gives a combinatorial interpretation of $N_S(m,n)$.

\begin{thm}\label{main-1}
For any integer  $m$ and any positive integer $n$,
 $N_S(m,n)$ equals the number of doubly marked partitions of $n$ with spt-crank $m$.
\end{thm}

 For example, for $n=4$, the sixteen $S$-partitions
of $4$, their spt-cranks  and the ten doubly marked partitions of $4$ and their spt-cranks  are listed in  Table 2.1. It can be checked that
\[N_S(3,4)=N_S(-3,4)=N_S(2,4)=N_S(-2,4)=1,\]
and
\[N_S(1,4)=N_S(-1,4)=N_S(0,4)=2.\]

\begin{table}[h]
\[\begin{array}{ccc|cc}
S\text{-partition}& \text{weight}&\text{spt-crank}&\text{doubly marked partition}&\text{spt-crank}\\[2pt]
((1),(1,1,1),\emptyset)&+1&3& ((1,1,1,1),1,1)&3\\[2pt]\hline
((1),(2,1),\emptyset)&+1&2&((2,1,1),1,1)&2\\[2pt]\hline
((1),(1,1),(1))&+1&1&((3,1),1,1)&1\\[2pt]
((1),(3),\emptyset)&+1&1&((2,2),1,2)&1\\[2pt]
((2,1),(1),\emptyset)&-1&1&&\\[2pt]
((2),(2),\emptyset)&+1&1&&\\[2pt]\hline
((1),(2),(1))&+1&0&((2,2),1,1)&0\\[2pt]
((1),(1),(2))&+1&0&((4),1,4)&0\\[2pt]
((3,1),\emptyset,\emptyset)&-1&0&&\\[2pt]
((4),\emptyset,\emptyset)&+1&0&&\\[2pt]\hline
((1),(1),(1,1))&+1&-1&((2,2),2,2)&-1\\[2pt]
((1),\emptyset,(3))&+1&-1&((4),1,3)&-1\\[2pt]
((2,1),\emptyset,(1))&-1&-1&&\\[2pt]
((2),\emptyset,(2))&+1&-1&&\\[2pt]\hline
((1),\emptyset,(2,1))&+1&-2&((4),1,2)&-2\\[2pt]\hline
((1),\emptyset,(1,1,1))&+1&-3&((4),1,1)&-3
\end{array}
\]
\caption{$S$-partitions and doubly marked partitions.}
\end{table}

The proof of  Theorem \ref{main-1} relies on the following generating function of $N_S(m,n)$ given by Andrews, Garvan and Liang \cite{Andrews-Garvan-Liang-2011}.
\begin{thm}
\begin{align}
&\sum_{m=-\infty}^{+\infty}\sum_{n\geq 0}N_S(m,n)z^mq^n\nonumber \\[3pt]
&\quad \quad =1+\sum_{m=0}^{+\infty}z^m
\sum_{j=0}^{+\infty}\frac{q^{j^2+mj+2j+m+1}}{(q;q)_{j+m}}\sum_{h=0}^{j}{j
\brack h}\frac{q^{h^2+h}}{(q;q)_h(1-q^{m+1+j+h})}\nonumber\\[3pt]
&\quad \quad \quad \quad +\sum_{m=1}^{+\infty}z^{-m}\sum_{j=m}^{+\infty}\frac{q^{j^2-mj+2j-m+1}}{(q;q)_{j-m}}\sum_{h=0}^{j}{j
\brack h}\frac{q^{h^2+h}}{(q;q)_h(1-q^{j-m+1+h})}.\label{spt-c}
\end{align}
\end{thm}

Let $Q_{m,n}$ denote  the set of doubly marked partitions of $n$
with spt-crank $m$. We aim to show  that  for $m\geq 0$
\begin{equation}\label{gf-q-p}
\sum_{n\geq 1}\sum_{(\lambda,s,t)\in Q_{m,n}}q^{|\lambda|}=\sum_{j=0}^{+\infty}\frac{q^{j^2+mj+2j+m+1}}{(q;q)_{j+m}}\sum_{h=0}^{j}{j
\brack h}\frac{q^{h^2+h}}{(q;q)_h(1-q^{m+1+j+h})},
\end{equation}
and
\begin{equation}\label{gf-q-n}
\sum_{n\geq 1}\sum_{(\lambda,s,t)\in Q_{-m,n}}q^{|\lambda|}=\sum_{j=m}^{+\infty}\frac{q^{j^2-mj+2j-m+1}}{(q;q)_{j-m}}
\sum_{h=0}^{j}{j
\brack h}\frac{q^{h^2+h}}{(q;q)_h(1-q^{j+h+1-m})}.
\end{equation}

To do this, we first represent a doubly marked partition $(\lambda,s,t)$ as a pair of partitions $(\alpha,\beta)$. Let $(\lambda,s,t)$ be a doubly marked partition of $n$ with spt-crank $m$, we define a pair of partitions $(\alpha,\beta)=\psi(\lambda,s,t)$ as follows.
  \begin{equation}\label{defi-alpha}
   \alpha=(\lambda_1-t+s-1,\ldots,\lambda_{\lambda_s'}-t+s-1,\lambda_{\lambda
  _s'+1},\ldots,\lambda_\ell)
  \end{equation}
  and
  \begin{equation}\label{defi-beta}
  \beta=(\lambda'_s,\lambda'_{s+1},\ldots, \lambda'_t)
   \end{equation}
  as illustrated in Figure 2.2. By the definition of a
   doubly marked partition, $\beta$ is a partition with equal parts. Let $P_{n}$ denote the set of pairs of partitions $(\alpha,\beta)$ of $n$ where $\beta$ is a partition with equal parts. We first claim that this representation is unique for the doubly marked partitions in $Q_{m,n}$. We then characterize the image set of $\psi$. Thus we can use the
   image set of $\psi$ to compute the generating function of
    $Q_{m,n}$.

   \input{lambdac.TpX}

  \begin{lem}\label{psi-i}
 Given integer $m$ and positive integer $n$, then the map $\psi$ is an injection from the set $Q_{m,n}$ to the set $P_{n}$. In other words, any doubly marked partition $(\lambda,s,t)$ in $Q_{m,n}$ can
 be uniquely represented by $(\alpha,\beta)$ in $P_{n}$.
 \end{lem}

 \pf Let $(\lambda,s,t)$ and $(\bar{\lambda},\bar{s},\bar{t})$ be two doubly marked partitions of $n$ with spt-crank $m$ such that
   \begin{equation}\label{invpsi-1}
   \psi(\lambda,s,t)=\psi(\bar{\lambda},\bar{s},\bar{t})= (\alpha,\beta).
   \end{equation}
   We proceed to show that $(\lambda,s,t)=(\bar{\lambda},\bar{s},\bar{t})$.

 By the construction of $\psi$ and relation (\ref{invpsi-1}),
 we have    \[\lambda=\bar{\lambda}=(\alpha_1+\ell(\beta), \alpha_2+\ell(\beta),\ldots, \alpha_{\beta_1}+\ell(\beta), \alpha_{\beta_1+1}, \ldots, \alpha_\ell).\]
 It remains to show that $s=\bar{s}$ and $t=\bar{t}$.

By the definition \eqref{def-sptc} of the spt-crank, we have \[g-\lambda_g+\ell(\beta)-1=m\]
 and
 \[\bar{g}-\bar{\lambda}_{\bar{g}}+\ell(\beta)-1=m,\]
 where
 \begin{equation}\label{defi-g}
 g=\lambda'_{s}-s+1
 \end{equation}
  and
 \begin{equation}\label{defi-g-bar}
 \bar{g}=\bar{\lambda}'_{\bar{s}}-\bar{s}+1.
 \end{equation}
It follows that
\begin{equation}\label{psi-inv-g}
 g-\lambda_g=\bar{g}-\bar{\lambda}_{\bar{g}}.
\end{equation}
We claim that  $g=\bar{g}$. Assume to the contrary that $g\not=\bar{g}$.
Without loss of generality, we may assume that $g>\bar{g}$. By \eqref{psi-inv-g}, we see that $\lambda_g>\bar{\lambda}_{\bar{g}}$. Since $\lambda=\bar{\lambda}$, we get
$\lambda_g>\lambda_{\bar{g}}$ which implies $g<\bar{g}$, contradicting
 the assumption that $g>\bar{g}$. So we have verified that $g=\bar{g}$.

Now that $g=\bar{g}$, by  \eqref{defi-g} and \eqref{defi-g-bar}, we see that
$\lambda'_{s}-s+1=\bar{\lambda}'_{\bar{s}}-\bar{s}+1$, that is,
 \[\lambda'_{s}-s=\bar{\lambda}'_{\bar{s}}-\bar{s}.\]
 Using the same argument as in the derivation  of $g=\bar{g}$
 from relation (\ref{psi-inv-g}),  we deduce that $s=\bar{s}$.  By the construction of $\psi$, we see that
 \[\ell(\beta)=t-s+1=\bar{t}-\bar{s}+1.\]
 But we have shown that $s=\bar{s}$, it follows that $t=\bar{t}$. Thus we
  reach the conclusion that $(\lambda,s,t)=(\bar{\lambda},\bar{s},\bar{t}).$ This completes the proof.\qed

 The following lemma will be used to  characterize the image set  of the map  $\psi$. To present this lemma,  we need to recall the definitions
of the rank-set  and the $m$-Durfee rectangle of a partition.  Let $\lambda=(\lambda_1,\lambda_2,\ldots, \lambda_\ell)$ be an ordinary partition. The rank-set  of $\lambda$ was introduced by Dyson \cite{Dyson-1989}
as an infinite sequence
 \[ [- \lambda_1, 1-\lambda_2, \ldots, j-\lambda_{j+1},\ldots,\ell-1-\lambda_{\ell}, \ell, \ell+1,\ldots] .\]
 For example, the rank-set of $\lambda=(5,5,4,3,1)$ is $ [-5,-4,-2,0,3,5,6,7,8,\ldots].$

When $m$ is an integer, the $m$-Durfee rectangle  of a partition $\lambda$ was introduced by Gordon and Houten \cite{Gordon-Houten-1968}, see also, Andrews \cite{Andrews-1971},
as the largest $(m+j)\times j$ rectangle   contained in the Ferrers diagram of $\lambda$, see Figure 2.3, where   $j$ is called the width of the $m$-Durfee rectangle  of $\lambda$.  An $m$-Durfee rectangle reduces to a Durfee square when $m=0$.  It should be noticed that $m$ may be negative.

\input{mdurfee.TpX}

It is worth mentioning that for   a partition $\lambda$ with  $\ell(\lambda)\leq m$,
 there is no $m$-Durfee rectangle. In this case,
 we adopt a convention that the  $m$-Durfee rectangle has
 no columns, that is, we set $j=0$.

\begin{lem}\label{prop-2-alpha}
Let $(\lambda,s,t)$ be a doubly marked partition of $n$ with spt-crank $m$,
 and let  $(\alpha,\beta)=\psi(\lambda,s,t)$. Then we have the following
 properties.
\begin{itemize}
\item [{\rm (1)}] The rank-set of $\alpha$ contains $m${\rm ;}
\item[{\rm (2)}] Let $j$ be the width of the $m$-Durfee rectangle of $\alpha$ and $h$ be the maximum integer such that $\alpha_{j+m+1+h}\geq h$, then
    \[\beta_1=j+m+1+h.\]
\end{itemize}
\end{lem}

\pf We first show that $m$ appears in the rank-set of $\alpha$.
Since $(\lambda, s,t)$ has spt-crank $m$, by definition,
we have
\begin{equation}\label{crank-dm}
g-\lambda_g+t-s=m,
\end{equation}
where $g=\lambda'_s-s+1$.
Since $s\geq 1$, we have $g\leq \lambda_s'$. Noting that $s\leq D(\lambda)$, we get $\lambda'_s\geq s$. It follows that $g\geq 1$.  So we have
$1\leq g\leq \lambda_s'$. By  the construction of $\psi$,    we see that
\begin{equation}\label{alpha-g-lambda-g}
\lambda_{g}-(t-s+1)=\alpha_{g}.
\end{equation}
Substituting  \eqref{alpha-g-lambda-g} into \eqref{crank-dm}, we obtain
  \begin{equation}\label{rs-alpha}
  g-1-\alpha_{g}=m,
  \end{equation}
  which implies $m$ appears the rank-set of $\alpha$.

  We continue to show that $\beta_1=j+m+h+1$.   Recall that $j$ is the width of  the $m$-Durfee rectangle of $\alpha$.
   Since   $m$ appears in the rank-set of $\alpha$, we find that
  \begin{equation}\label{rs-jm}
  j+m-\alpha_{j+m+1}=m.
  \end{equation}
 From \eqref{rs-alpha} and \eqref{rs-jm}, we deduce that
    \begin{equation}\label{defi-t-1}
    g=j+m+1.
    \end{equation}
Since $g=\lambda'_s-s+1$, it follows from  \eqref{defi-t-1}  that
 \begin{equation}\label{temp-1-1}
 \lambda_s'-s=j+m.
 \end{equation}
 From the construction of $\psi$, we see that $\beta_1=\lambda_s'$.  Hence
  \[\beta_1=s+j+m.\]

  We claim that $s=h+1$. By the choice of $h$,  it suffices  to show that
   $\alpha_{j+m+s}\geq s-1$ and ${\alpha}_{j+m+s+1}<s.$
 By \eqref{temp-1-1}, we have
  \[\alpha_{j+m+1+s}=\alpha_{\lambda_s'+1},\]
  and
  \[\alpha_{j+m+s}=\alpha_{\lambda_s'}.\]
   From the construction of $\psi$, we see that \[\alpha_{\lambda_s'+1}=\lambda_{\lambda
  '_s+1},\]
  and
  \[\alpha_{\lambda'_s}=\lambda_{\lambda'_s}-t+s-1.\]
  Examining the Ferrers diagram of $\lambda$, we see that $\lambda_{\lambda'_s+1}<s$. Consequently,
  \begin{equation}\label{claim-2-1}
  \alpha_{j+m+s+1}=\alpha_{\lambda_s'+1}=\lambda_{\lambda
  '_s+1}<s.
  \end{equation}
  Since $\lambda_t'=\lambda_s'$, it can be seen
   that $\lambda_{\lambda_s'}\geq t$. So we deduce that
  \begin{equation}\label{claim-2-2}
  \alpha_{j+m+s}=\alpha_{\lambda'_s}=\lambda_{\lambda'_s}-t+s-1\geq t-t+s-1=s-1.
  \end{equation}
 Combining \eqref{claim-2-1} and \eqref{claim-2-2}, we deduce that $h=s-1$.
 So we conclude that $\beta_1=j+m+h+1$.  This completes the proof. \qed

  It turns out that the properties in Lemma \ref{prop-2-alpha} are sufficient
  to characterize the image of the map $\psi$.

  \begin{thm}\label{lemdmpp} Given an integer $m$ and a positive integer $n$, let $V_{m,n}$ denote the set of pairs of partitions $(\alpha,\beta)$ of $n$ satisfying the conditions in Lemma \ref{prop-2-alpha}. The map $\psi$ is a bijection between the set $Q_{m,n}$ and the set $V_{m,n}$.
  \end{thm}

\pf By Lemma \ref{psi-i} and Lemma \ref{prop-2-alpha}, it suffices to construct the inverse map $\varphi$ such that for all   $(\alpha,\beta)$ in $V_{m,n}$, we have $\psi(\varphi(\alpha,\beta))=(\alpha,\beta).$

Given a pair of partitions $(\alpha,\beta)$ in $V_{m,n}$,  we   construct a doubly marked partition $(\lambda,s,t)$. Recall that  $h$ is the maximum integer such that $\alpha_{j+m+1+h}\geq h$. Let $s=h+1$, $t=h+\ell(\beta)$ and
 $\lambda=(\alpha_1+\ell(\beta), \alpha_2+\ell(\beta),\ldots, \alpha_{\beta_1}+\ell(\beta), \alpha_{\beta_1+1}, \ldots, \alpha_\ell).$
 Then it can be seen that  $(\lambda,s,t)$ is a doubly marked partition. Set
 $\varphi(\alpha,\beta)=(\lambda,s,t)$. By the definitions of $\psi$ and $\varphi$, it can be verified
  that  $\psi(\varphi(\alpha,\beta))=(\alpha,\beta)$ for all $(\alpha,\beta)$ in $V_{m,n}$. The detailed steps are omitted. This completes the proof. \qed

We are now in a position to complete the proof of Theorem \ref{main-1}.

\noindent{\it Proof of Theorem \ref{main-1}.}  By the bijection $\psi$ in Theorem \ref{lemdmpp} between the set $Q_{m,n}$ and the set $V_{m,n}$, we see that for all integers $m$ and  positive integers $n$,
 \[\sum_{(\lambda,s,t)\in Q_{m,n}}q^{|\lambda|}=\sum_{(\alpha,\beta)\in V_{m,n}}q^{|\alpha|+|\beta|}.\]
To prove (\ref{gf-q-p}) and (\ref{gf-q-n}), it suffices to show that for $m\geq 0$, we have
\begin{equation}\label{gf-pp-p}
\sum_{n\geq 1}\sum_{(\alpha,\beta)\in V_{m,n}}q^{|\alpha|+|\beta|}=\sum_{j=0}^{+\infty}\frac{q^{j^2+mj+2j+m+1}}{(q;q)_{j+m}}\sum_{h=0}^{j}{j
\brack h}\frac{q^{h^2+h}}{(q;q)_h(1-q^{m+1+j+h})},
\end{equation}
and
\begin{equation}\label{gf-pp-n}
\sum_{n\geq 1}\sum_{(\alpha,\beta)\in V_{-m,n}}q^{|\alpha|+|\beta|}=\sum_{j=m}^{+\infty}\frac{q^{j^2-mj+2j-m+1}}{(q;q)_{j-m}}
\sum_{h=0}^{j}{j
\brack h}\frac{q^{h^2+h}}{(q;q)_h(1-q^{j+h+1-m})}.
\end{equation}
We first consider  \eqref{gf-pp-p}. Let $V_{m,n}^{j,h}$ denote the set of pairs of partitions $(\alpha,\beta)$ in $V_{m,n}$ such that the width of the $m$-Durfee rectangle of $\alpha$ is $j$ and $h$ is the maximum number
satisfying $\alpha_{j+m+1+h}\geq h$. Note that $\alpha_{j+m+1}=j$.
 Hence we have $0\leq h\leq j$. Now, we have the following
 relation for the sum on the left hand side
 of (\ref{gf-pp-p})
\begin{equation}\label{gf-pp-p-2}
\sum_{n\geq 1}\sum_{(\alpha,\beta)\in V_{m,n}}q^{|\alpha|+|\beta|}=\sum_{j=0}^{+\infty}\sum_{h=0}^j \sum_{n\geq 1}\sum_{(\alpha,\beta)\in V_{m,n}^{j,h}}q^{|\alpha|+|\beta|}.
\end{equation}
We shall show that
\begin{equation}\label{gf-pp-p-1}
\sum_{n\geq 1}\sum_{(\alpha,\beta)\in V_{m,n}^{j,h}}q^{|\alpha|+|\beta|}=\frac{q^{j^2+mj+2j+m+1}}{(q;q)_{j+m}}{j
\brack h}\frac{q^{h^2+h}}{(q;q)_h(1-q^{m+1+j+h})}.
\end{equation}
For any $(\alpha,\beta)$ in  $V_{m,n}^{j,h}$, by Lemma \ref{prop-2-alpha}, we see that $\beta$ is a partition with each part equal to $j+m+1+h$. So   the generating function for $\beta$   is given by
\begin{equation}\label{gf-beta}
\frac{q^{h+j+m+1}}{1-q^{m+1+j+h}}.
\end{equation}
To derive the generating function of $\alpha$, we  decompose the Ferrers
diagram of  $\alpha$ into four regions  as shown in Figure 2.4.
The generating function for $\alpha$ can be determined
by computing the generating function of each region.

\input{alphac.TpX}

Region $A$ forms to a partition into $j+m$ parts for which each part equals $j$.
The generating function for this region is $q^{(m+j)\times j}$.
Region $B$ is a partition with only one part $j$, whose generating function is $q^j$.

Region $C$ forms to a partition with each part not exceeding $j$.
 Let $\gamma$ be the partition in this region. Since $h$ is the maximum number  such that $\alpha_{j+m+1+h}\geq h$, we see that  the size of the Durfee square of $\gamma$ is  $h$. Therefore,  we may divide $\gamma$ into
 three partitions $\gamma^1$, $\gamma^2$ and $\gamma^3$:
\begin{itemize}
\item[(1)] $\gamma^1$ is the Durfee square of $\gamma$, which is of size
$h\times h$.
\item[(2)] $\gamma^2$ is the partition formed by the parts to the
 right of $\gamma^1$, which is a partition into at most $h$ parts, each part not exceeding  $j-h$;
\item[(3)] $\gamma^3$ is the partition formed by the parts below $\gamma^1$, which is a partition with each  part not exceeding $h$.
\end{itemize}
Hence the generating function for all possible partitions $\gamma$ in
this region is given by
\[  q^{h^2}  {j
\brack h}\frac{1}{(q;q)_h}, \]
where
\[{n+m \brack m}=\frac{(q;q)_{n+m}}{(q;q)_n(q;q)_m}\]
is the generating function for partitions with at most $m$ parts, each part not exceeding $n$, see \cite[p.35]{Andrews-1976}.

Region $D$ forms a partition into at most $m+j$ parts. So the
 generating function for possible partitions in this region equals
 \[\frac{1}{(q;q)_{m+j}}.\]

Taking all the regions into consideration, we obtain the generating function of all possible partitions $\alpha$ satisfying the
constraints of the set $V_{m,n}^{j,h}$
\begin{align}\label{gf-alpha}
 q^{(m+j)\times j}\cdot q^j \cdot  q^{h^2}  {j
\brack h}\frac{1}{(q;q)_h} \cdot \frac{1}{(q;q)_{m+j}}.
\end{align}

It should be noticed that for given $m,n,j$ and $h$,
to form a pair $(\alpha, \beta)$ in $V_{m,n}^{j,h}$, the
choices  $\alpha$ and $\beta$ are independent. Hence the
 generating function for $(\alpha,\beta)$ in $ V_{m,n}^{j,h}$ is
 the product of the generating functions of $\alpha$ and $\beta$ subject
 to their individual constraints. Thus   \eqref{gf-pp-p-1} follows from \eqref{gf-beta} and \eqref{gf-alpha}. So we arrive at  \eqref{gf-pp-p}.

  We now turn to the proof of (\ref{gf-pp-n}). The generating function for pairs of partitions in $V_{-m,n}$ can be
   computed in the same way as the derivation of (\ref{gf-pp-p}).
  Let $V_{-m,n}^{j,h}$ denote the set of pairs of partitions $(\alpha,\beta)$ in $V_{-m,n}$ such that the width of the $-m$-Durfee rectangle of $\alpha$ is $j$, or the size of the $-m$-Durfee rectangle is $(j-m)\times j$. Assume
  that $h$ is the maximum number  such that $\alpha_{j+m+1+h}\geq h$.
  Since $j\geq m$, we have
\begin{equation}\label{gf-pp-n-2}
\sum_{n\geq 1}\sum_{(\alpha,\beta)\in V_{-m,n}}q^{|\alpha|+|\beta|}=\sum_{j=m}^{+\infty}\sum_{h=0}^j \sum_{n\geq 1}\sum_{(\alpha,\beta)\in V_{-m,n}^{j,h}}q^{|\alpha|+|\beta|}.
\end{equation}
Using the same argument as in the proof of \eqref{gf-pp-p-1},  we can deduce that
\begin{equation}\label{gf-pp-n-1}
\sum_{n\geq 1}\sum_{(\alpha,\beta)\in V_{-m,n}^{j,h}}q^{|\alpha|+|\beta|}=\frac{q^{j^2-mj+2j-m+1}}{(q;q)_{j-m}}{j
\brack h}\frac{q^{h^2+h}}{(q;q)_h(1-q^{j-m+1+h})}.
\end{equation}
Combining \eqref{gf-pp-n-2} and \eqref{gf-pp-n-1}, we obtain \eqref{gf-pp-n}.    This completes the proof.
 \qed

\section{Marked partition and doubly marked partitions}

In this section, we establish a correspondence between marked partitions and doubly marked partitions, so that one can divide the set of marked partitions of $5n+4$ and $7n+5$
into five and seven equinumerous classes by employing the spt-crank
of doubly marked partitions.

\begin{thm}\label{main-2} There is a bijection $\Delta$ between the set of marked partitions $(\mu,k )$ of $n$ and the set of doubly marked partitions $(\lambda, s, t)$ of $n$.
\end{thm}

To prove the above theorem, we need to use the notation $(\lambda, s,t)$
to mean a partition $\lambda$ with two distinguished columns in the Ferrers diagram. In other words, we no longer assume that $(\lambda, s,t)$
is a doubly marked partition unless it is explicitly stated.
Let $Q_n$ denote the set of doubly marked partitions of $n$,
 and let
 \[U_n=\{(\lambda,s,t) \ |\    |\lambda|=n,\, 1\leq s\leq D(\lambda),\, 1\leq t\leq \lambda_1\}.\]
 Obviously, $Q_n\subseteq U_n$.

Before we give a description of the bijection $\Delta$, we  introduce a transformation $\tau$ from   $U_n\setminus Q_n$ to $U_n$.

\noindent{\it The transformation $\tau$:} Assume that $(\lambda,s,t) \in U_n\setminus Q_n$, that is, $\lambda$ is an ordinary partition of $n$ with
two distinguished columns $s$ and $t$ such that $1\leq s\leq D(\lambda)$  and either $1\leq t<s$ or $\lambda'_{s}>\lambda'_t.$  We wish to construct a partition $\mu$ with two distinguished columns $a$ and $b$.  Let  $p$ be the maximum integer such that $\lambda_p'=\lambda_s'$. Define
\begin{equation}\label{phi-gamma}
\delta=(\lambda_1-p+s-1,\lambda_2-p+s-1,\ldots, \lambda_{\lambda_s'}-p+s-1,\lambda_{\lambda_s'+1},\ldots, \lambda_\ell).
\end{equation}
 Set $a$ to be   the minimum integer such that $\delta_a<\lambda'_s$ and
\begin{equation}\label{def-mu}
\mu=(\delta_1,\ldots,\delta_{a-1},\lambda'_s,\ldots, \lambda'_p, \delta_{a},\ldots, \delta_\ell).
\end{equation}
If $t<s$, then set $b=t$ and if $\lambda_s'>\lambda_t'$, then   set $b=t-p+s-1$.  Define $\tau(\lambda,s,t)=(\mu,a,b)$. Figure 3.1 gives an illustration of the map $\tau$.
\input{phi.TpX}

The conditions in the following lemma will be used to  characterize the image set of the map $\tau$.

\begin{lem}\label{lem-1} Assume that $(\lambda,s,t)\in U_n\setminus Q_n$,
and denote  $\tau(\lambda,s,t)$ by $(\mu, a,b)$.
  We have  $(\mu, a,b) \in U_n$.  Furthermore, if  $a=1$, then we have $\mu'_b>s(\mu').$
\end{lem}
\pf  We first show that   $(\mu,a,b) \in U_n$, that is,  we need to verify that $1\leq a\leq D(\mu)$ and $1\leq b\leq \mu_1$.

It is clear that $a\geq 1$. We proceed to show that $a \leq D(\mu)$. To this end, we first prove that $\delta_{\lambda_s'}<\lambda_{s}'$. Then we show that $a \leq D(\mu)$ can be deduced from the fact that $\delta_{\lambda_s'}<\lambda_{s}'$.

By the definition of $\delta$ in \eqref{phi-gamma}, we see that \begin{equation}\label{del-ine}
\delta_{\lambda_s'}=\lambda_{\lambda_s'}-p+s-1<\lambda_{\lambda_s'}.
\end{equation}
Since $1\leq s\leq D(\lambda)$, that is, $\lambda_s'\geq D(\lambda)$,
we deduce that
\begin{equation}\label{lam-ine}
\lambda_{\lambda_s'}\leq \lambda_s'.
\end{equation}
Combining \eqref{del-ine} and \eqref{lam-ine} yields $\delta_{\lambda_s'}<\lambda_s'$.

Recall that $a$ is the minimum integer such that $\delta_a<\lambda_{s}'$.
But we have shown that  $\delta_{\lambda_s'}<\lambda_{s}'$, this implies
that $ a\leq \lambda_{s}'$. On the other hand,
 by the construction of $\tau$, we find that $ \lambda_s'=\mu_a$.
 So we deduce that $a\leq \mu_a$, that is,  $a\leq D(\mu)$.
 This completes the proof of the assertion that $1\leq a \
\leq D(\mu)$.

Next, we continue to prove that $1\leq b\leq \mu_1$. There are two cases.

\noindent
Case 1:  $1\leq t<s$. By the construction of $\tau$, we have $b=t$ and
\begin{equation}\label{temp-1}
 \mu_1\geq \lambda_s'.
\end{equation}
Since $s\leq D(\lambda)$, we get
\begin{equation}\label{temp-2}
\lambda_s'\geq s>t.
\end{equation}
Combining \eqref{temp-1} and \eqref{temp-2}, we deduce that $t< \mu_1$. Since $b=t$ and $t\geq 1$, we conclude that $1\leq b<\mu_1$.

\noindent
Case 2:    $\lambda_s'>\lambda_t'$. By the construction of $\tau$, we have $b=t-p+s-1$ and
\begin{equation}\label{phi-lem-temp1}
\mu_1\geq \delta_1=\lambda_1-p+s-1.
\end{equation}
Since  $\lambda_p'=\lambda_s'$ and $\lambda_s'>\lambda_t'$, we have $t>p$,
and so $b=t-p+s-1\geq s\geq 1$.
Using \eqref{phi-lem-temp1} and the fact that $t\leq \lambda_1,$  we obtain $\mu_1\geq t-p+s-1=b$. It follows that  $1\leq b\leq \mu_1$.

Up to now, we have shown  that $(\mu,a,b) \in U_n.$
Finally, we prove that if $a=1$, then $\mu'_b>s(\mu').$
 We now assume that $a=1$, and we claim that in this case
 \begin{equation}\label{temp-3}
b \leq \lambda_1-p+s-1.
\end{equation}

By the choice of $a$, if $a=1$, then we have
\begin{equation}\label{temp-det}
\delta_1<\lambda'_s
\end{equation} and
 \begin{equation}\label{temp-5}
 \mu=(\lambda'_s,\ldots, \lambda'_p, \delta_{1},\ldots, \delta_\ell).
 \end{equation}

 To prove the claim, we consider the following two cases.

\noindent Case 1:  $1\leq t<s$.  By the construction of $\tau$, we have $b=t \leq s-1$.  Since $p\leq \lambda_1$, we see that
 $$ b=t\leq s-1\leq \lambda_1-p+s-1.$$

 \noindent Case 2:  $\lambda_s'>\lambda_t'$. By the construction of $\tau$, we find that $b=t-p+s-1$. Using the fact that $t\leq \lambda_1$,   we  get
 $$b=t-p+s-1\leq \lambda_1-p+s-1.$$
So the claim is proved.

Combining \eqref{temp-det} and \eqref{temp-5}, we get $\delta_1<\mu_1$,
or equivalently, $\mu'_{\delta_1}>\mu_{\mu_1}'$. Note that $\mu_{\mu_1}'=s(\mu')$, so we have $\mu'_{\delta_1}>s(\mu')$. On the other hand, from the definition \eqref{phi-gamma} of $\delta$, we have $\delta_1=\lambda_1-p+s-1$. By the claim that $b \leq \lambda_1-p+s-1$,
we obtain that  $b\leq \delta_1$, this yields  $\mu'_b\geq \mu'_{\delta_1}$.
So we reach the conclusion that   $\mu'_b>s(\mu')$.  This completes the proof.  \qed

The following theorem gives the image set $W_n$ of the transformation $\tau$,
and it shows that $\tau$ is bijection between $U_n\setminus Q_n$ and
 $W_n$.

\begin{thm}\label{lem-2}Given a positive integer $n$, let
 \begin{equation}\label{defi-wn}
 W_n=\{(\mu,a,b)\  | \ (\mu,a,b)\in U_n\ \text{and} \  \mu_b'>s(\mu')  \text{ whenever } a=1\}.
 \end{equation}
Then the transformation $\tau$ is a bijection between   $U_n\setminus Q_n$ and $W_n$.
\end{thm}

\pf By Lemma \ref{lem-1},  it suffices to construct a map  $\sigma$ defined on $W_n$ such that for all $(\lambda,s,t) \in U_n\setminus Q_n$, we have $\sigma(\tau(\lambda,s,t))=(\lambda,s,t)$ and for all $(\mu,a,b) \in W_n$, we have $\tau(\sigma(\mu,a,b))=(\mu,a,b).$

Let $(\mu,a,b)\in W_n$, we wish to construct a partition $\lambda$ with two distinguished columns $s$ and $t$. Let $r$ be the maximum integer such that $\mu_r=\mu_a$. Define
\[\gamma=(\mu_1,\ldots, \mu_{a-1},\mu_{r+1},\ldots,\mu_\ell).\]
  Set  $s$ to be the minimum integer such that $\gamma_s'<\mu_a$, and
\[\lambda=(\gamma_1+r-a+1,\ldots, \gamma_{\mu_a}+r-a+1,\gamma_{\mu_a+1},\ldots, \gamma_{\ell}).\]
If $b<s$, then we set $t=b$. Otherwise, we set $t=b+r-a+1$. Define $\sigma(\mu,a,b)=(\lambda,s,t)$. Using the same argument as in the proof of Lemma \ref{lem-1}, we deduce that $\sigma(\mu,a,b)\in U_n\setminus Q_n$.

By the constructions of $\tau$ and $\sigma$, it is straightforward to check that $\sigma(\tau(\lambda,s,t))=(\lambda,s,t)$ for all $(\lambda,s,t) \in U_n\setminus Q_n$ and $\tau(\sigma(\mu,a,b))=(\mu,a,b)$ for all $(\mu,a,b) \in W_n$. The details  are omitted.  This completes the proof. \qed

We now describe the bijection $\Delta$ in Theorem \ref{main-2} based on the bijection $\tau$.

\noindent{\it The definition of $\Delta \colon$} Let $(\mu, k)$ be a marked partition of $n$, we wish to construct a doubly marked partition $(\lambda,s,t)$ of $n$.

 We first consider  $(\mu', 1, k)$.
 If $(\mu', 1, k)$ is already a doubly marked partition, then
 there is nothing to be done and we just set $(\lambda,s,t)=(\mu',1,k)$. Otherwise, we iteratively apply the map $\tau$ to $(\mu',1,k)$
 until we get a doubly marked partition $(\lambda,s,t)$. We shall show that
 this process terminates and it is reversible.

 For example, let $n=6$, $\mu=(2,1,1,1,1)$ and $k=5$. We have $\mu'=(5,1)$.
 Note that $(\mu', 1, k)=((5,1),1,5)$, which  is not a doubly marked partition.
 It can be checked that  $\tau(\mu', 1,k)=
       ((4,2),2,4)$, which is not a doubly marked partition.
 Repeating this process, we get $\tau((4,2),2,4)=((3,2,1),2,3)$,
  and $\tau ((3,2,1),2,3)=((2,2,1,1),2,2)$, which is eventually  a doubly marked partition.  See Figure 3.2. Thus, we obtain
 \[\Delta((2,1,1,1,1),5)=((2,2,1,1),2,2).\]

 \input{triexa.TpX}

The following lemma shows that the map $\Delta$ is well-defined.

\begin{lem} \label{trian-lem}
The map $\Delta$ is well-defined, that is, for each marked partition $(\mu,k)$, there exists $i$ such that $\tau^i(\mu',1,k)$ is a doubly marked partition.
\end{lem}

\pf Assume to the contrary that  there exists a marked partition $(\mu,k)$  of
$n$ such that for any $i\geq 0$, $\tau^i(\mu',1,k)$ is not a doubly marked partition of $n$.
Let $(\lambda^{(i)},s^{(i)},t^{(i)})=\tau^i(\mu',1,k)$. By Lemma \ref{lem-1}, we see that $\lambda^{(i)}$ is an ordinary partition of $n$, $s^{(i)}$ and $t^{(i)}$ are both bounded by $n$. Thus the set
$$\{(\lambda^{(i)},s^{(i)},t^{(i)})| \ i \geq 0\}$$
is finite. So there exist integers $\ell$ and $m$ such that $\ell<m$ and $
(\lambda^{(\ell)},s^{(\ell)},t^{(\ell)})=(\lambda^{(m)},s^{(m)},t^{(m)}),
$
that is,
\begin{equation}\label{initial}
\tau^{\ell}(\mu',1,k)=\tau^m(\mu',1,k).
\end{equation}
We may choose  $\ell$ to the minimum integer such that $\tau^\ell(\mu',1,k)=\tau^m(\mu',1,k)$ for some $m>\ell$.   We claim that $\ell\geq 1$, that is, there does not exist $m\geq 1$ such that
\begin{equation}\label{triangle-tem}
(\mu',1,k)=\tau^m(\mu',1,k).
\end{equation}
 Denote $(\mu',1,k)$ by $(\lambda,a,b)$, so that we have $\lambda'=\mu$. Since $(\mu, k)$ is a marked partition, that is, $\mu_k=s(\mu)$, we see that $\lambda_b'=\mu_k=s(\lambda')$.
 Since $a=1$, by the definition \eqref{defi-wn} of $W_n$, we see that $(\lambda,a,b)$ is   in $W_n$ if any only if $\lambda'_b>s(\lambda')$. So we deduce that
 $(\mu', 1,k)$ is not in $W_n$. On the other hand,  by Theorem \ref{lem-2}, we see that $W_n$ is the image-set of $\tau$. Since $\tau^m(\mu',1,k)$ lies  in the image-set of $\tau$, it follows that $(\mu',1,k)\neq \tau^m(\mu',1,k)$ for any $m\geq 1$.   This proves that $\ell \geq 1$.

By the choice of $\ell$, we see that  for any $m>\ell$,    \begin{equation}\label{tau-l-m}
\tau^{\ell-1}(\mu',1,k)\neq \tau^{m-1}(\mu',1,k).
\end{equation}
  By the assumption that $\tau^i(\mu',1,k)$ is not a doubly marked partition for any $i\geq 0$,  $\tau^{\ell-1}(\mu',1,k)$ and $\tau^{m-1}(\mu',1,k)$ are not doubly marked partitions.
   Since $\tau$ is a bijection, we obtain that
\[\tau^{\ell}(\mu',1,k)\neq \tau^{m}(\mu',1,k),\]
for any $m> \ell$,
 contradicting the choice of $\ell$. Hence we conclude that
 $\Delta$ is well-defined. This completes the proof.   \qed

We are now ready to complete the proof of Theorem \ref{main-2}.

\noindent{\it Proof of Theorem \ref{main-2}:}
We have given the description of the map $\Delta$ and have shown that
$\Delta$ is well-defined. It remains to show that
$\Delta$ is reversible. To this end, we construct a map $\Lambda$ defined on the set of doubly marked partitions of $n$ and we shall show that it is
 the inverse map of $\Delta$.

 The map $\Lambda$ can be described as follows.
 Let $(\lambda,s,t)$ be a doubly marked partition of $n$. We aim to
 construct a marked partition $(\mu,k)$.
  If $s=1$  and $\lambda_t'=s(\lambda')$, we set  $(\mu,k)= (\lambda',t)$.
  If $s>1$ or $\lambda_t' > s(\lambda')$, we can iteratively use the inverse map $\tau^{-1}$ to transform $(\lambda,s,t)$
 into an ordinary partition $\delta$ with two distinguished columns $a$ and $b$ such that $a=1$ and $\delta_b'=s(\delta')$. Set $(\mu,k)=(\delta',b)$.
  Finally, define $\Lambda(\lambda,s,t)=(\mu,k).$

Parallel to the proof of the fact that $\Delta$ is well-defined, it can be shown that the map $\Lambda$ is well-defined. The details are omitted. Furthermore, since $\tau$ is a bijection, it is routine to check that
 $\Lambda$ is the inverse map of $\Delta$. This completes the proof.\qed

Employing the bijection $\Delta$ and the spt-crank for doubly
marked partitions, one can  divide the set of marked partitions of $5n+4$ and $7n+5$ into five and seven equinumerous classes.

 For example, for $n=4$, we have $spt(4)=10$. The ten marked partitions
of $4$, the corresponding doubly marked partitions, and the spt-crank modulo $5$ are listed in  Table 3.1.

{\small \begin{table}[h]\label{table-1}
 \[\begin{array}{c|c|c|c}
(\mu,k)& (\lambda,s,t)=\Delta(\mu,k) & c(\lambda,s,t) &c(\lambda,s,t) \mod 5\\[2pt] \hline
((4),1) & ((1,1,1,1),1,1)& 3& 3\\[2pt]
((3,1),2)& ((3,1),1,1)& 1& 1\\[2pt]
((2,2),1)& ((2,2),1,1)& 0& 0\\[2pt]
((2,2),2)& ((2,2),1,2)& 1& 1\\[2pt]
 ((2,1,1),2)& ((2,1,1),1,1)& 2 &2\\[2pt]
((2,1,1),3)& ((2,2),2,2)& -1 &4\\[2pt]
((1,1,1,1),1)& ((4),1,1) & -3 &2\\[2pt]
((1,1,1,1),2)& ((4),1,2) & -2 &3\\[2pt]
((1,1,1,1),3)& ((4),1,3) & -1 &4\\[2pt]
((1,1,1,1),4)& ((4),1,4) & 0 &0
\end{array}\]

 \caption{ The case for $n=4$.}
\end{table}}

\newpage

For $n=5$, we have $spt(5)=14$. The fourteen marked partitions
of $5$, the corresponding doubly marked partitions, and the spt-crank modulo $7$ are listed in  Table 3.2.

{\small
\begin{table}[h]\label{table-2}
\[\begin{array}{c|c|c|c}
(\mu,k)& (\lambda,s,t)=\Delta(\mu,k) & c(\lambda,s,t) & c(\lambda,s,t)\mod 7\\[2pt] \hline
((5),1) & ((1,1,1,1,1),1,1)& 4& 4\\[2pt]
((4,1),2)& ((4,1),1,1)& 1& 1\\[2pt]
((3,2),2)& ((3,1,1),1,1)& 2& 2\\[2pt]
((3,1,1),2)& ((3,2),1,1)& 0& 0\\[2pt]
((3,1,1),3)& ((3,2),1,2)& 1& 1\\[2pt]
((2,2,1),3)& ((2,2,1),1,1)& 2& 2\\[2pt]
((2,1,1,1),2)& ((2,1,1,1),1,1)& 3 &3\\[2pt]
((2,1,1,1),3)& ((3,2),2,2)& -2 &5\\[2pt]
((2,1,1,1),4)& ((2,2,1),2,2)& -1 &6\\[2pt]
((1,1,1,1,1),1)& ((5),1,1) & -4 &3\\[2pt]
((1,1,1,1,1),2)& ((5),1,2) & -3 &4\\[2pt]
((1,1,1,1,1),3)& ((5),1,3) & -2 &5\\[2pt]
((1,1,1,1,1),4)& ((5),1,4) & -1 &6\\[2pt]
((1,1,1,1,1),5)& ((5),1,5) & 0 &0
\end{array}\]

 \caption{ The case for $n=5$.}
\end{table}
}

\noindent{\bf Acknowledgments.} This work was supported by the 973 Project, the PCSIRT Project of the Ministry of Education and the National Science Foundation of China.

\end{document}